\documentclass[conference, letterpaper, 10pt]{ieeeconf}
\IEEEoverridecommandlockouts                              

\usepackage{xspace,amssymb,epsfig,subfigure,syntonly}
\usepackage{epsfig,amsmath,color}
\usepackage{xspace,syntonly}
\usepackage{url}
\usepackage{algorithm,algorithmic}
\usepackage{mathtools}
\usepackage{cite}
\usepackage{pstool}

\newcommand{\utwi}[1]{\mbox{\boldmath $#1$}}

\newcommand{\diag}{{\textrm{diag}}}

\newcommand{\cD}{{\cal D}}

\newcommand{\cL}{{\cal{L}}}
\newcommand{\cN}{{\cal N}}

\newcommand{\cG}{{\cal G}}

\newcommand{\cC}{{\cal C}}
\newcommand{\cE}{{\cal E}}

\newcommand{\cF}{{\cal F}}

\newcommand{\cM}{{\cal M}}

\newcommand{\cY}{{\cal Y}}

\newcommand{\bc}{{\bf c}}
\newcommand{\ba}{{\bf a}}
\newcommand{\bb}{{\bf b}}

\newcommand{\be}{{\bf e}}
\newcommand{\bbf}{{\bf f}}
\newcommand{\bg}{{\bf g}}

\newcommand{\bm}{{\bf m}}

\newcommand{\bp}{{\bf p}}
\newcommand{\bq}{{\bf q}}
\newcommand{\br}{{\bf r}}

\newcommand{\bx}{{\bf x}}
\newcommand{\bu}{{\bf u}}

\newcommand{\bv}{{\bf v}}
\newcommand{\bi}{{\bf i}}
\newcommand{\bz}{{\bf z}}
\newcommand{\by}{{\bf y}}

\newcommand{\bB}{{\bf B}}

\newcommand{\bJ}{{\bf J}}
\newcommand{\bH}{{\bf H}}

\newcommand{\bR}{{\bf R}}

\newcommand{\bX}{{\bf X}}

\newcommand{\bY}{{\bf Y}}

\newcommand{\bgamma}{{\utwi{\gamma}}}
\newcommand{\blambda}{{\utwi{\lambda}}}

\newcommand{\brho}{{\utwi{\rho}}}

\newcommand{\bPhi}{{\utwi{\Phi}}}

\newcommand{\bxi}{{\utwi{\xi}}}

\newcommand{\bmu}{{\utwi{\mu}}}
\newcommand{\bzeta}{{\utwi{\zeta}}}

\newcommand{\sfH}{\textsf{H}}
\newcommand{\sfT}{\textsf{T}}

\usepackage{epstopdf}

\begin{document}

\newtheorem{definition}{Definition}
\newtheorem{remark}{Remark}
\newtheorem{proposition}{Proposition}
\newtheorem{lemma}{Lemma}
\newtheorem{theorem}{Theorem}
\def\HS{\hspace{\fontdimen2\font}}
\font\myfont=cmr12 at 16pt

\IEEEoverridecommandlockouts
\overrideIEEEmargins


\title{\LARGE \bf
Design of Distributed Controllers Seeking Optimal Power Flow Solutions Under Communication Constraints
}

\author{Emiliano Dall'Anese, Andrea Simonetto, and Sairaj Dhople 
\thanks{E. Dall'Anese is  with the National Renewable Energy Laboratory (NREL), Golden, CO, USA. A Simonetto is with the Universit\'e Catholique de Louvain, Louvain-la-Neuve, Belgium. S. Dhople is 
with the Department of Electrical and Computer Engineering, University of Minnesota, Minneapolis, USA. E-mail: emiliano.dallanese@nrel.gov, andrea.simonetto@uclouvain.be, sdhople@umn.edu.}
\thanks{The work of E. Dall'Anese was supported in part by the Laboratory Directed Research and Development Program at NREL. S. Dhople was supported in part by the National Science Foundation under the CAREER award 1453921 and CyberSEES grant 1442686.}
}

\maketitle

\maketitle

\begin{abstract}
This paper focuses on power distribution networks featuring distributed energy resources (DERs), and develops controllers that drive the DER output powers to solutions of time-varying AC optimal power flow (OPF) problems. The design of the controllers is grounded on primal-dual-type methods for regularized Lagrangian functions, as well as linear approximations of the AC power-flow equations. Convergence and OPF-solution-tracking capabilities are established while acknowledging: i) communication-packet losses, and ii) partial updates of control signals. The latter case is particularly relevant since it enables an asynchronous operation of the controllers where the DER setpoints are updated at a fast time scale based on local voltage measurements, and information on the network state is utilized if and when available, based on communication constraints. As an application, the paper considers distribution systems with a high penetration level of photovoltaic systems, and demonstrates that the proposed framework provides fast voltage-regulation capabilities, while enabling the near real-time pursuit of AC OPF solutions.
\end{abstract}

\section{Introduction}
\label{sec:Introduction}

Centralized and distributed AC optimal power flow (OPF) approaches have been developed  for distribution systems to compute optimal  setpoints for distributed energy resources (DERs), so that power losses and voltage deviations are minimized and economic benefits to utility and end-users are maximized. It is well-known that the AC OPF is a \emph{nonconvex} (and, in fact, NP-hard) nonlinear program. Centralized approaches utilize off-the-shelf solvers for nonlinear programs~\cite{Paudyal11}, or, leverage convex relaxation and approximation techniques to obtain convex surrogates~\cite{Farivar12,OID,Robbins15,swaroop2015linear}. Distributed solution methods tap into the decomposability of the Lagrangian function associated with convex surrogates of the OPF, and utilize iterative primal-dual-type methods to decompose the solution of the OPF task across DERs, utility, and possibly aggregators~\cite{Dallanese-TSG13,Erseghe14,Robbins15}. Either way, in the presence of (fast) changing load, ambient, and network conditions~\cite{Bank13}, traditional centralized and distributed OPF schemes may  offer decision making capabilities that do not match the  dynamics of distribution systems. Particularly, during the time required to collect data from all the nodes of the network (e.g., loads),  solve the OPF, and subsequently dispatch the setpoints, the load, ambient, and network conditions may have already changed. In this case, the DER output powers would be consistently regulated around outdated setpoints, leading to suboptimal system operation and violation of relevant electrical limits. This motivates the development of online OPF strategies that leverage the opportunities for fast-feedback offered by power-electronics-interfaced DERs to enable the near real-time pursuit of solutions of AC OPF problems,  while ensuring adaptability to fast-changing conditions~\cite{OPFpursuit,LowOnlineOPF, DhopleDKKT15}.  

Prior efforts in this direction include the continuous-time feedback controllers that seek Karush-Kuhn-Tucker conditions for economic dispatch optimality  for bulk systems developed in~\cite{Jokic_JEPES}.  Modified automatic generation and frequency control methods that incorporate optimization objectives corresponding to DC OPF problems have been proposed for bulk power systems in, e.g.,~\cite{NaLi_ACC14}. Focusing on AC OPF models, online solution approaches include  the heuristic based on saddle-point-flow method utilized in~\cite{Elia-Allerton13}, the online OPF proposed in~\cite{LowOnlineOPF} for distribution systems with a tree topology, and the distributed dual (sub)-gradient scheme developed in~\cite{DhopleDKKT15} for (un)balanced distribution systems. Overall, the convergence results in~\cite{Jokic_JEPES,LowOnlineOPF,DhopleDKKT15} hinge on a time scale separation where cost and constraints of the  OPF problem change slower than the controller dynamics.  A centralized controller is developed in~\cite{AndreayOnlineOpt} based on  gradient algorithms; it is shown that the DER setpoints convergence on average to an optimal solution. 

A distributed control architecture that enables DERs to track the solution of fast-changing OPF solutions is developed in~\cite{OPFpursuit}. Stability and tracking capabilities are  characterized in terms of bounds between the DER output powers and the optimal trajectory set forth by the \emph{time-varying} OPF problem. The present paper significantly broadens the approach~\cite{OPFpursuit} by considering realistic scenarios where communication constraints lead to \emph{asynchronous} and \emph{partial} updates of the control signals. Similar to~\cite{OPFpursuit},  control synthesis is based on suitable linear approximations of the AC power-flow equations as well as Lagrangian regularization methods. OPF-solution tracking  is  established for the cases where: i) communication-packet losses lead to asynchronous updates of the control signals; and ii) DER setpoints are updated at a fast time scale based on local voltage measurements, and information on state of the the remaining part of the network is utilized if and when available, based on communication constraints. This setup allows controllers to ensure that OPF constraints are met, while relaxing the requirements on the supporting communication infrastructure.

\section{Time-varying optimal power flow}
\label{sec:preliminariesandsystemmodel}

Consider a distribution feeder\footnote{Upper-case (lower-case) boldface letters will be used for matrices (column vectors); $(\cdot)^\sfT$ for transposition; $(\cdot)^*$
  complex-conjugate; and, $(\cdot)^\sfH$ complex-conjugate transposition; $\Re\{\cdot\}$ and $\Im\{\cdot\}$ denote the real and
  imaginary parts of a complex number, respectively; $\mathrm{j} := \sqrt{-1}$ the imaginary unit; and $|\cdot|$ denotes the absolute value of a number or the cardinality of a set. For $x \in \mathbb{R}$, function $[x]_+$ is defined as $[x]_+ := \max\{0,x\}$. For a given $N \times 1$ vector $\bx \in \mathbb{R}^N$, $\|\bx\|_2 := \sqrt{\bx^\sfH \bx}$; and, $\diag(\bx)$ returns a $N \times N$ matrix with the elements of $\bx$ in its diagonal. Further,  $\mathrm{proj}_{\cY}\{\bx\}$ denotes the projection of $\bx$ onto the convex set $\cY$. Given a given matrix $\bX \in \mathbb{R}^{N\times M}$, $x_{m,n}$  denotes its $(m,n)$-th entry. $\nabla_{\bx} f(\bx)$ returns the gradient vector of $f(\bx)$ with respect to $\bx \in \mathbb{R}^N$. Finally, $\mathbf{1}_N$ denotes the $N \times 1$ vector with all ones, and $\mathbf{0}_N$ denotes the $N \times 1$ vector with all zeros.} comprising $N+1$ nodes collected in the set $\cN \cup \{0\}$, $\cN := \{1,\ldots,N\}$, and lines represented by the set of
edges $\cE := \{(m,n)\} \subset \cN  \cup \{0\} \times \cN  \cup \{0\}$. Assume that  the temporal domain is discretized as $t = k \tau$, where $k \in \mathbb{N}$ and $\tau > 0$ is small enough to capture fast variations of loads and ambient conditions. Let $V_n^k \in \mathbb{C}$ and $I_n^k \in \mathbb{C}$ denote the phasors for the line-to-ground voltage and the current injected at node $n$ over the $k$th instant, respectively, and define the $N$-dimensional complex vectors  $\bv^k := [V_1^k, \ldots, V_N^k]^\sfT \in \mathbb{C}^{N}$ and $\bi^k := [I_1^k, \ldots, I_N^k]^\sfT \in
\mathbb{C}^{N}$. Node $0$ denotes the distribution transformer, and it is taken to be the slack bus. Using Ohm's and Kirchhoff's circuit laws, it follows that $\bi^k = V_0^k \overline \by^k + \bY^k \bv^k$, where $\bY^k \in \mathbb{C}^{N \times N}$ and $\overline \by^k \in \mathbb{C}^{N \times 1}$ are formed based on the  network topology and the $\pi$-equivalent circuit of the lines (see e.g.,~\cite{kerstingbook}). 

Inverter-interfaced DERs such as photovoltaic (PV) systems, small-scale wind turbines, and energy storage systems are assumed to be located at nodes $\cG \subseteq \cN$, $N_\cG := |\cG|$. The real and reactive powers at the AC side of inverter $i \in \cG$ at each time $k \tau$ are denoted as $P_{i}^k$ and $Q_i^k$, respectively, and are confined within the DER operating region $(P_{i}^k, Q_{i}^k) \in \cY_i^k$. The set $\cY_i^k$ captures hardware as well as operational constraints, and is assumed to be convex and compact. For example, for PV inverters, this set is given by $\cY_i^k =  \{({P}_{i}^k, {Q}_{i}^k ) \hspace{-.1cm} :  {P}_{i}^{\textrm{min}}  \leq {P}_{i}^k  \leq  P_{\textrm{av},i}^k, ({Q}_{i}^k)^2  \leq  S_{i}^2 - ({P}_{i}^k)^2\}$, where $P_{\textrm{av},i}^{k}$ denotes the real power available at time $k$ and $S_{i}$ is the capacity of the inverter~\cite{Farivar12,OID}.  In commercial-scale HVAC systems, set $\cY_i^k $ captures the power consumption range of a fan (i.e., a variable frequency drive). For future developments, let $\bu_i^k := [P_i^k, Q_i^k]^\sfT$ collect the real and reactive setpoints for DER $i$ at time $k$, and define the set $\cY^k := \cY_1^k \times \ldots \cY_{N_\cG}^k$. Finally, for each node $i$, let $P_{\ell,i}^k$ and $Q_{\ell,i}^k$ denote the inflexible real and reactive power demand, respectively, at time $k$. 

\subsection{Time-varying Target Optimization Problem}
\label{sec:systemmodel}

To bypass challenges related to nonconvexity and NP-hardness of the OPF task, and facilitate the design of low-complexity controllers that afford implementation on microcontrollers that accompany power-electronics interfaces of inverters, the present paper leverages suitable linear approximations of the AC power-flow equations. To this end, collect the voltage magnitudes $\{|V_i^k|\}_{i \in \cN}$ in the vector $\brho^k :=  [|V_1^k|, \ldots, |V_N^k|]^\sfT \in \mathbb{R}^{N}$. Then, given pertinent matrices $\bR^k, \bB^k, \bH^k, \bJ^k \in \mathbb{R}^{N \times N}$ and vectors $\bb^k, \ba^k \in \mathbb{C}^N$, one can obtain approximate power-flow relations whereby voltages  are \emph{linearly} related to the injected real and reactive powers as 
\begin{subequations} 
\label{eq:approximateVoltages}
\begin{align} 
\bv^k & \approx \bH^k \bp^k + \bJ^k \bq^k + \bb^k \label{eq:approximateV} \\
\brho^k & \approx \bR^k \bp^k + \bB^k \bq^k + \ba^k, \label{eq:approximate} 
\end{align}
\end{subequations}
where $p^k_n = P_n^k - P_{\ell,n}^k$, $q^k_n = Q_n^k- Q_{\ell,n}^k $ if $n \in \cG$ and   $p^k_n = - P_{\ell,n}^k$, $q^k_n = - Q_{\ell,n}^k $ if $n \in \cN \backslash \cG$. Matrices $\bR^k, \bB^k, \bH^k, \bJ^k \in \mathbb{R}^{N \times N}$ and vectors $\bb^k, \ba^k \in \mathbb{C}^N$ can be obtained as described   in e.g.,~\cite{sairaj2015linear,swaroop2015linear,bolognani2015linear}, and can be time-varying to reflect e.g., changes in the topology and voltage linearization points. It is worth pointing out that through~\eqref{eq:approximateV}--\eqref{eq:approximate}  approximate linear relationships for power losses and power flows as a function of $\{P_i^k, Q_i^k \}_{i \in \cG}$ can be readily derived~\cite{sairaj2015linear,swaroop2015linear}. 

Denote as $V^{\mathrm{min}}$ and $V^{\mathrm{max}}$ minimum and maximum, respectively, voltage service limits, and let the cost $\sum_{n \in \cG} f_n^k(\bu_n^k)$ capture possibly time-varying DER-oriented objectives (e.g., cost of/reward for ancillary service provisioning~\cite{Farivar12,OID}, or feed-in tariffs), and/or system-level performance metrics (e.g., power losses and/or deviations from the nominal voltage profile~\cite{OID}). With these definitions, and based on~\eqref{eq:approximateV}--\eqref{eq:approximate}, an approximate convex AC OPF problem can be formulated as~\cite{swaroop2015linear,OPFpursuit}:
\begin{subequations} 
\label{Pmg2}
\begin{align} 
 \mathrm{(P1}^k \mathrm{)}  \hspace{1.8cm} & \hspace{-1.7cm} \min_{\{\bu_i \}_{i \in \cG} } \,\, \sum_{i \in \cG} f_i^k(\bu_i) \label{mg-cost2} \\ 
& \hspace{-1.5cm} \mathrm{subject\,to}   \nonumber  \\ 
& \hspace{-0.8cm} g^k_n(\{\bu_i\}_{i \in \cG}) \leq 0 , \hspace{1.45cm} \forall n \in \cM \label{mg-volt1} \\
& \hspace{-0.8cm}  \bar{g}^k_n(\{\bu_i\}_{i \in \cG}) \leq 0 ,  \hspace{1.45cm} \forall n \in \cM \label{mg-volt2} \\
& \hspace{-0.8cm} \bu_i \in  \cY_i^k ,  \hspace{2.75cm} \forall \, i \in \cG\ , \label{mg-PVp2} 
\end{align}
\end{subequations}
where $\cM \subseteq \cN$ is a set of nodes selected to enforce voltage regulation throughout the feeder, $M := |\cM|$, and
\begin{subequations} 
\label{eq:v} 
\begin{align} 
g^k_n(\{\bu_i\}_{i \in \cG}) & := V^{\mathrm{min}} - c_n^k \nonumber \\
&   - \sum_{i \in \cG} [r_{n,i}^k (P_i - P_{\ell,i}^k) + b_{n,i}^k (Q_i - Q_{\ell,i}^k)] \label{eq:vmin} \\
\bar{g}^k_n(\{\bu_i\}_{i \in \cG}) &:= \sum_{i \in \cG} [r_{n,i}^k (P_i - P_{\ell,i}^k) + b_{n,i}^k (Q_i - Q_{\ell,i}^k)] \nonumber \\
& + c_n^k  - V^{\mathrm{max}} \, , \label{eq:vmax} 
\end{align}
\end{subequations} 
with $c_n^k := a_n^k - \sum_{i \in \cN \backslash \cG} (r_{n,i}^k P_{\ell,i}^k + b_{n,i}^k Q_{\ell,i}^k)$. Regarding~\eqref{Pmg2}, the following  assumptions are made.  

\vspace{.1cm}

\noindent \emph{Assumption~1}.  Functions $f_i^k(\bu_i)$ are convex and continuously differentiable for each $i \in \cG$ and $k \geq 0$. Define further the gradient map  $\bbf^k(\bu) := [\nabla_{\bu_1}^\sfT f_1^k(\bu_1), \ldots, \nabla_{\bu_{N_\cG}}^\sfT f_{N_\cG}^k(\bu_{N_\cG})]^\sfT$.
Then, it is assumed that $\bbf^k: \mathbb{R}^{2 N_\cG} \rightarrow \mathbb{R}^{2 N_\cG}$ is Lipschitz continuous
with constant $L$ over $\cY^k$ for all $k \geq 0$. \hfill $\Box$ 

\vspace{.1cm}

\noindent \emph{Assumption~2}. For all $k \geq 0$, there exist a set of feasible power injections $\{\hat{\bu}_i\}_{i \in \cG} \in \cY^{k}$ such that $g^k_n(\{\hat{\bu}_i\}_{i \in \cG}) \leq 0$ and $\bar{g}^k_n(\{\hat{\bu}_i\}_{i \in \cG}) \leq 0$, for all $n \in \cM$.  \hfill $\Box$ 

\vspace{.1cm}

Regarding \emph{Assumption~2}, notice that functions $g^k_n(\{\hat{\bu}_i\}_{i \in \cG})$ and $\bar{g}^k_n(\{\hat{\bu}_i\}_{i \in \cG})$ are linear [cf.~\eqref{eq:v}]; hence, Slater's condition does not require strict inequalities~\cite{BoVa04}. From the compactness of set $\cY^k$, and under \emph{Assumptions~1} and \emph{2}, problem~\eqref{Pmg2} is convex and strong duality holds. Further, there exists an optimizer $\{\bu_i^{\mathrm{opt},k}\}_{i \in \cG}$, $\forall \,\, k \geq 0$. For future developments, let $\bg^k(\bu) \in \mathbb{R}^M$ and $\bar{\bg}^k(\bu) \in \mathbb{R}^M$ be a vector stacking all functions $g^k_n(\{\bu_i\}_{i \in \cG})$ and $\bar{g}^k_n(\{\bu_i\}_{i \in \cG})$.

\subsection{Objective}
\label{sec:objective}

Problem $ \mathrm{(P1}^k \mathrm{)} $ represents a convex approximation of the AC OPF task. Constraints~\eqref{mg-volt1}--\eqref{mg-volt2} are utilized to enforce voltage regulation, while~\eqref{mg-PVp2} models DER hardware constraints. It is worth pointing out that the problem $\mathrm{(P1}^k \mathrm{)}$ specifies OPF targets that corresponds to a specific time instant $k \tau$; accordingly, in the presence of (fast) changing load, ambient, and network conditions, repeated solutions of $\mathrm{(P1}^k \mathrm{)}$ for $k \in \mathbb{N}$ would ideally produce optimal reference setpoint \emph{trajectories} for the DER $\{\bu_n^{\mathrm{opt},k}, k \in \mathbb{N}\}$. However, traditional centralized and distributed solution approaches may not be able to collect network data (e.g., loads),  solve $ \mathrm{(P1}^k \mathrm{)} $, and subsequently dispatch setpoints within $\tau$ seconds, and may consistently regulate the power-outputs $\{P_i^k, Q_i^k \}_{i \in \cG}$ around outdated setpoints. This motivates the development of  controllers that continuously regulate the DER output powers around points that one would have  if $ \mathrm{(P1}^k \mathrm{)} $ could be solved instantaneously.  

Particularly, let $\by^k = \cF(\{\bu_i^k\}_{i \in \cG})$ represent an AC power-flow solution for given DER output powers $\{\bu_i^k \}_{n \in \cG}$, with vector $\by^k$ collecting relevant electrical quantities such as voltages and power flows (averaged over one AC cycle)~\cite{Dorfler14,DhopleDKKT15,LowOnlineOPF,kerstingbook}. Further, let $\cC_i\{\cdot, \by^k\}$ describe an update rule for the setpoints of DER $i$. Then, given the following closed loop-system
\begin{subequations}
\label{eq:Objectivecontroller}
\begin{align}
\bu_i^k & = \cC_i(\bu_i^{k-1}, \by^k) , \hspace{.5cm} \forall i \in \cG  \label{eq:ObjectivecontrollerC} \\
\by^k & = \cF(\{\bu_i^k\}_{i \in \cG})
\end{align}
\end{subequations}
the goal is to design the controllers $\{\cC_i\{\cdot, \cdot\}\}_{i \in \cG}$ so that the DER output powers $\{\bu_i^k\}_{i \in \cG}$ are driven to the solution $\{\bu_i^{\mathrm{opt},k}\}_{i \in \cG}$ of the time-varying OPF problem $ \mathrm{(P1}^k \mathrm{)} $. 

\section{Design of feedback controllers}
\label{sec:opfPursuit}

\subsection{Preliminaries}
\label{sec:preliminaries}

The synthesis of the controllers leverages primal-dual methods applied to regularized Lagrangian functions~\cite{Koshal11,SimonettoGlobalsip2014}. To this end,  let $\bgamma := [\gamma_1, \ldots, \gamma_M]^\sfT$ and $\bmu := [\mu_1, \ldots, \mu_M]^\sfT$ collect the Lagrange multipliers associated with~\eqref{mg-volt1} and~\eqref{mg-volt2}, respectively, and consider the following augmented Lagrangian function associated with $ \mathrm{(P1}^k \mathrm{)}$:  
\begin{align} 
& \cL_{\nu, \epsilon}^k(\bu^k, \bgamma, \bmu) := \sum_{i \in \cG} f_i^k(P_i, Q_i) + (P_i - P_{\ell,i}^k) (\check{\br}_i^k)^\sfT (\bmu - \bgamma) \nonumber \\
& + (Q_i - Q_{\ell,i}^k) (\check{\bb}_i^k)^\sfT (\bmu - \bgamma) + \bc^\sfT(\bmu - \bgamma) + \bgamma^\sfT \mathbf{1}_m V^{\mathrm{min}}  \nonumber\\ 
& - \bmu^\sfT \mathbf{1}_m V^{\mathrm{max}} 
 + \frac{\nu}{2} \sum_{n \in \cG} \|\bu_n^k\|_2^2 - \frac{\epsilon}{2} (\|\bgamma\|_2^2 + \|\bmu\|_2^2) 
\label{eq:lagrangianR}
\end{align}
where  $\check{\br}_i^k := [\{r_{j,i}^k\}_{j \in \cM}]^\sfT$ and $\check{\bb}_i^k := [\{b_{j,i}^k\}_{j \in \cM}]^\sfT$ are $M \times 1$ vectors collecting the entries of  $\bR^k$ and $\bB^k$ in their $i$th column and rows corresponding to nodes in $\cM$,  $\bc^k := [\{c_j^k\}_{j \in \cM}]^\sfT$, and  constants $\nu > 0$ and $\epsilon > 0$ appearing in the Tikhonov regularization terms are design parameters. Function~\eqref{eq:lagrangianR} is strictly convex  in the primal variables $\bu^{k} :=[ \bu_1^{k}, \ldots,  \bu_{N_\cG}^{k}]^\sfT$ and strictly concave in the dual variables $\bgamma, \bmu$. The upshot of~\eqref{eq:lagrangianR} is that gradient-based approaches can be applied to find an approximate solution to $\mathrm{(P1}^k \mathrm{)} $ with improved convergence properties~\cite{Koshal11,SimonettoGlobalsip2014}. Further, it allows one to drop the strict convexity assumption on the cost function $\{f_i^k(\bu_i)\}_{i \in \cG}$~\cite{NaLi_ACC14, Elia-Allerton13,DhopleDKKT15} and to avoid averaging primal and dual variables~\cite{OzdaglarSaddlePoint09}. Accordingly, consider the following saddle-point problem:
\begin{align} 
\label{eq:saddlepoint}
\max_{\blambda \in \mathbb{R}^{M}_+, \bmu \in \mathbb{R}^{M}_+} \min_{\bu \in \cY^k}  \quad \cL_{\nu, \epsilon}^k(\bu^k, \bgamma, \bmu) 
\end{align}
and denote as $\bu^{*,k} :=[ \bu_1^{*,k}, \ldots,  \bu_{N_\cG}^{*,k}]^\sfT, \bgamma^{*,k}, \bmu^{*,k}$ the \emph{unique} primal-dual optimizer of~\eqref{eq:lagrangianR}.  In general, the solutions of~\eqref{Pmg2} and the regularized saddle-point problem~\eqref{eq:saddlepoint} are expected to be different; however, the discrepancy between $\bu_i^{\textrm{opt},k}$ and $\bu_i^{*,k}$ can be bounded as in~\cite[Lemma~3.2]{Koshal11}, whereas bounds of the constraint violation are substantiated in~\cite[Lemma~3.3]{Koshal11}. These bounds are proportional to $\sqrt{\epsilon}$; therefore, the smaller $\epsilon$, the smaller is the discrepancy between $\bu_n^{\textrm{opt},k}$ and $\bu_n^{*,k}$. 

To track the time-varying optimizers $\bz^{*,k} := [(\bu^{*,k})^\sfT, (\bgamma^{*,k})^\sfT, (\bmu^{*,k})^\sfT]^\sfT$ of~\eqref{eq:saddlepoint},  consider  the following online primal-dual gradient method~\cite{SimonettoGlobalsip2014}:  
\begin{subequations} 
\label{eq:updateopt}
\begin{align} 
\bu_i^{k + 1} & =  \mathrm{proj}_{\cY_i}\left\{\bu_i^{k}  - \alpha \nabla_{\bu_i}  \cL_{\nu, \epsilon}^k(\bu, \bgamma, \bmu) |_{\bu_i^{k}, \bgamma^k, \bmu^k} \right\} \label{eq:updatev} \\
\gamma_n^{k + 1} & =   \mathrm{proj}_{\cD_\gamma} \left\{\gamma_n^{k} +  \alpha (g^k_n(\{\bu_i^k\}_{i \in \cG}) - \epsilon \gamma_n^{k}) \right\}  \label{eq:updategamma} \\
\mu_n^{k + 1} & =   \mathrm{proj}_{\cD_\mu} \left\{\mu_n^{k} +  \alpha (\bar{g}^k_n(\bu_i^k\}_{i \in \cG}) - \epsilon \mu_n^{k}) \right\}\ ,  \label{eq:updatemu}
\end{align}
\end{subequations} 
where $\alpha > 0$ is the  stepsize, and $\cD_\gamma, \cD_\mu \subset \mathbb{R}^+$ are compact convex sets that can be chosen as explained in~\cite{Koshal11}. Step~\eqref{eq:updatev} is computed for each $i \in \cG$, whereas~\eqref{eq:updategamma}--\eqref{eq:updatemu} are performed for each node $n \in \cM$. Convergence of the  iterates $\bz^{k} := [(\bu^{k})^\sfT, (\bgamma^{k})^\sfT, (\bmu^{k})^\sfT]^\sfT$ to $\bz^{*,k}$ is established in~\cite[Theorem~1]{SimonettoGlobalsip2014}, and hinges on the following assumptions related to the temporal variability of~\eqref{eq:saddlepoint}. 

\vspace{.1cm}

\noindent \emph{Assumption~3}. There exists a constant $\sigma_\bu \geq 0$ such that $\|\bu^{*,k+1} - \bu^{*,k}\| \leq \sigma_\bu$ for all $k \geq 0$. \hfill $\Box$ 

\vspace{.1cm}

\noindent \emph{Assumption~4}. There exist constants $\sigma_d \geq 0$ and $\sigma_{\bar{d}} \geq 0$ such that $|g^{k+1}_n(\bu^{*,k+1}) - g^{k}_n(\bu^{*,k}) | \leq \sigma_d$ and $|\bar{g}^{k+1}_n(\bu^{*,k+1}) - \bar{g}^{k}_n(\bu^{*,k}) | \leq \sigma_{\bar{d}}$, respectively, for all $n \in \cN$ and $k \geq 0$. \hfill $\Box$ 

\vspace{.1cm}

It can be shown that the conditions of \emph{Assumption~4} translate into bounds for the discrepancy between the optimal dual variables over two consecutive time instants; that is, $\|\bgamma^{*,k+1} - \bgamma^{*,k}\| \leq \sigma_{\bgamma}$ and $\|\bmu^{*,k+1} - \bmu^{*,k}\| \leq \sigma_{\bmu}$ with $\sigma_{\bgamma}$ and $\sigma_{\bmu}$ given by~\cite[Prop.~1]{SimonettoGlobalsip2014}. Overall, upon defining $\bz^{*,k}:=[(\bu^{*,k})^\sfT, (\bgamma^{*,k})^\sfT, (\bmu^{*,k})^\sfT]^\sfT$  it also follows that $\|\bz^{*,k+1} - \bz^{*,k}\| \leq \sigma_{\bz}$ for a given $\sigma_{\bz} \geq 0$. Under \emph{Assumptions~1--4}, convergence of~\eqref{eq:updateopt} are investigated in~\cite[Theorem~1]{SimonettoGlobalsip2014}.

Similar to~\cite{OPFpursuit}, in the next section the updates~\eqref{eq:updateopt} are modified to accommodate actionable feedback from the distribution system. The proposed framework broadens the approach of~\cite{OPFpursuit} by considering a more realistic scenario where communication constraints lead to \emph{asynchronous} and \emph{partial} updates of primal/dual variables.

Before proceeding, it is worth pointing out the following two facts: i) given that $\bg^k(\bu)$ and $\bar{\bg}^k(\bu)$  are linear in $\bu$ and $\cY^k$ is compact, it follows that there exists a constant $G$ such that $\|\nabla_\bu \bg^k(\bu)\|_2 \leq G$ and $\|\nabla_\bu \bar{\bg}^k(\bu)\|_2 \leq G$ for all $k \geq 0$. One can show that there exist constants $K > 0$ and $\bar{K}$ such that $\|\bg^k(\bu)\|_2 \leq K$ and $\|\bar{\bg}^k(\bu)\|_2 \leq \bar{K}$. Further, notice that $\|\bgamma^k_i\|_2 \leq D_\gamma$ and $\|\bmu^k_i\|_2 \leq D_\mu$ for given $D_\gamma, D_\mu >0$ by construction [cf.~\eqref{eq:updategamma}--\eqref{eq:updatemu}], and define the time-varying mapping $\bPhi^k$ as
\begin{equation*}
\label{phi_mapping}
 \bPhi^k: \{\bu^{k}, \bgamma^{k}, \bmu^{k}\} \mapsto 
 \left[\begin{array}{c}
 \nabla_{\bu_1}  \cL_{\nu, \epsilon}^k(\bu, \bgamma, \bmu) |_{\bu_1^{k}, \bgamma^k, \bmu^k} \\
 \vdots \\
 \nabla_{\bu_{N_{\cG}}}  \cL_{\nu, \epsilon}^k(\bu, \bgamma, \bmu) |_{\bu_{N_{\cG}}^{k}, \bgamma^k, \bmu^k}\\
- (g^k_1(\bu^k) - \epsilon \gamma_1^{k})  \\
\vdots \\
- (g^k_M(\bu^k) - \epsilon \gamma_M^{k})  \\
 - (\bar{g}^k_1(\bu^k) - \epsilon \mu_1^{k}) \\
 \vdots\\
  - (\bar{g}^k_M(\bu^k) - \epsilon \mu_M^{k})  \end{array}\right] . \hspace{-.1cm}
\end{equation*}
Then, the following holds.

\vspace{.1cm}

\begin{lemma}
\label{lemma-Phi}
The map $\bPhi^k$ is strongly monotone with constant $\eta = \min\{\nu,\epsilon\}$, and Lipschitz over $\cY^k \times \cD_\gamma  \times \cD_\mu$ with constant $L_{\nu, \epsilon} = \sqrt{(L+\nu+2G)^2 + 2(G+\epsilon)^2}$. \hfill $\Box$
\end{lemma}

\subsection{Feedback Controllers}
\label{sec:communicationConstraints}

With regards to the distributed optimization scheme~\eqref{eq:updateopt}, it is worth pointing out that: (i) functions $\{g^k_n(\bu^k)\}_{n \in \cM}$ and $\{\bar{g}^k_n(\bu^k)\}_{n \in \cM}$ capture the distance of the voltage magnitudes from the limits $V^{\mathrm{min}}$ and $V^{\mathrm{max}}$, respectively, of given setpoints $\bu^k$; (ii) to evaluate $g^k_n(\bu^k), \bar{g}^k_n(\bu^k)$ at the current points $\bu^k$ it is necessary to collect all loads  across the network  [cf.~\eqref{eq:v}]; and, (iii) all dual variables need to be collected at each DER $i \in \cG$ in order to carry out step~\eqref{eq:updatev}. 

\begin{figure}[t] 
\subfigure[]{\hspace{-.2cm}\includegraphics[width=9.0cm]{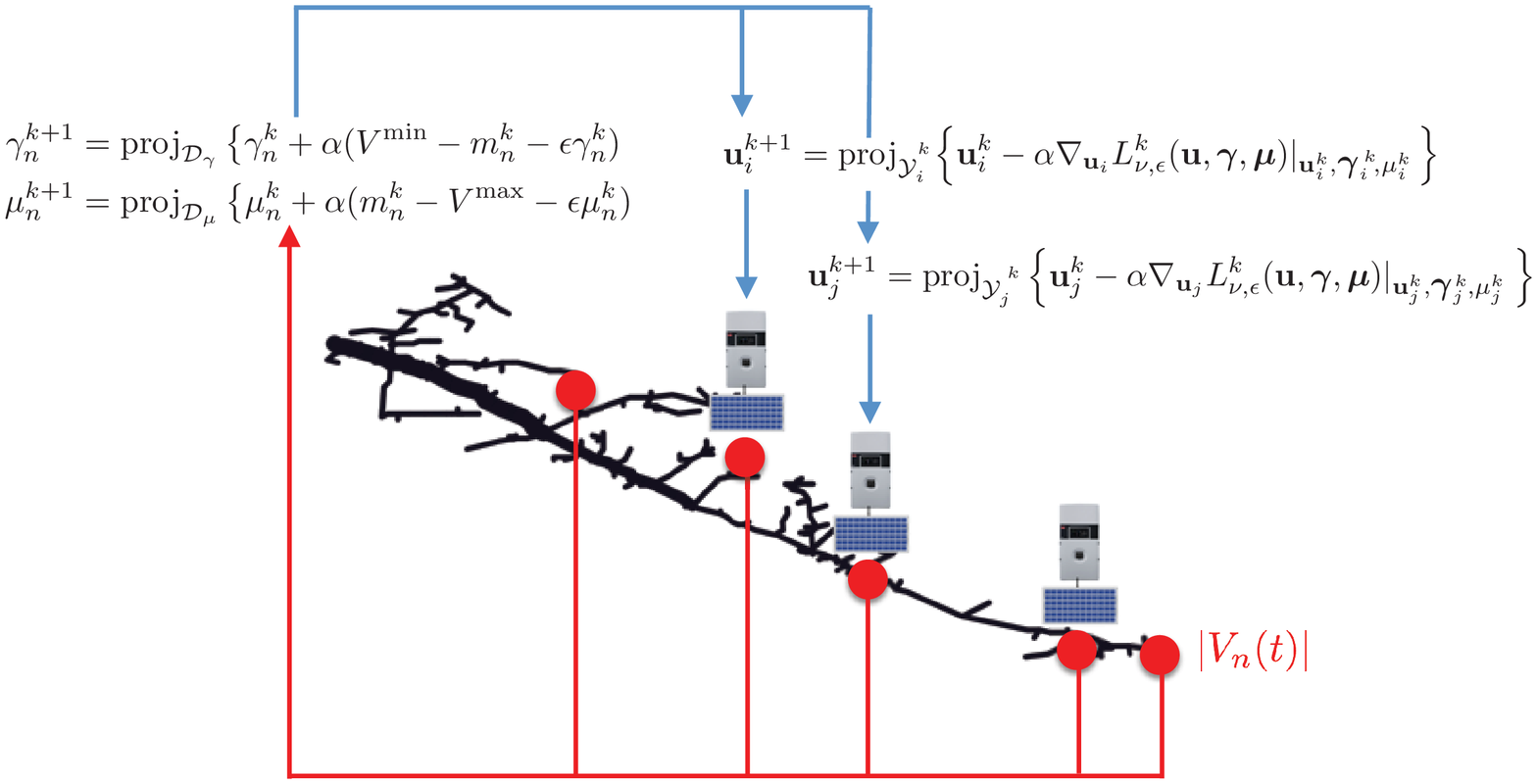} }\vspace{-.2cm}
\subfigure[]{\hspace{-.2cm}\includegraphics[width=9.0cm]{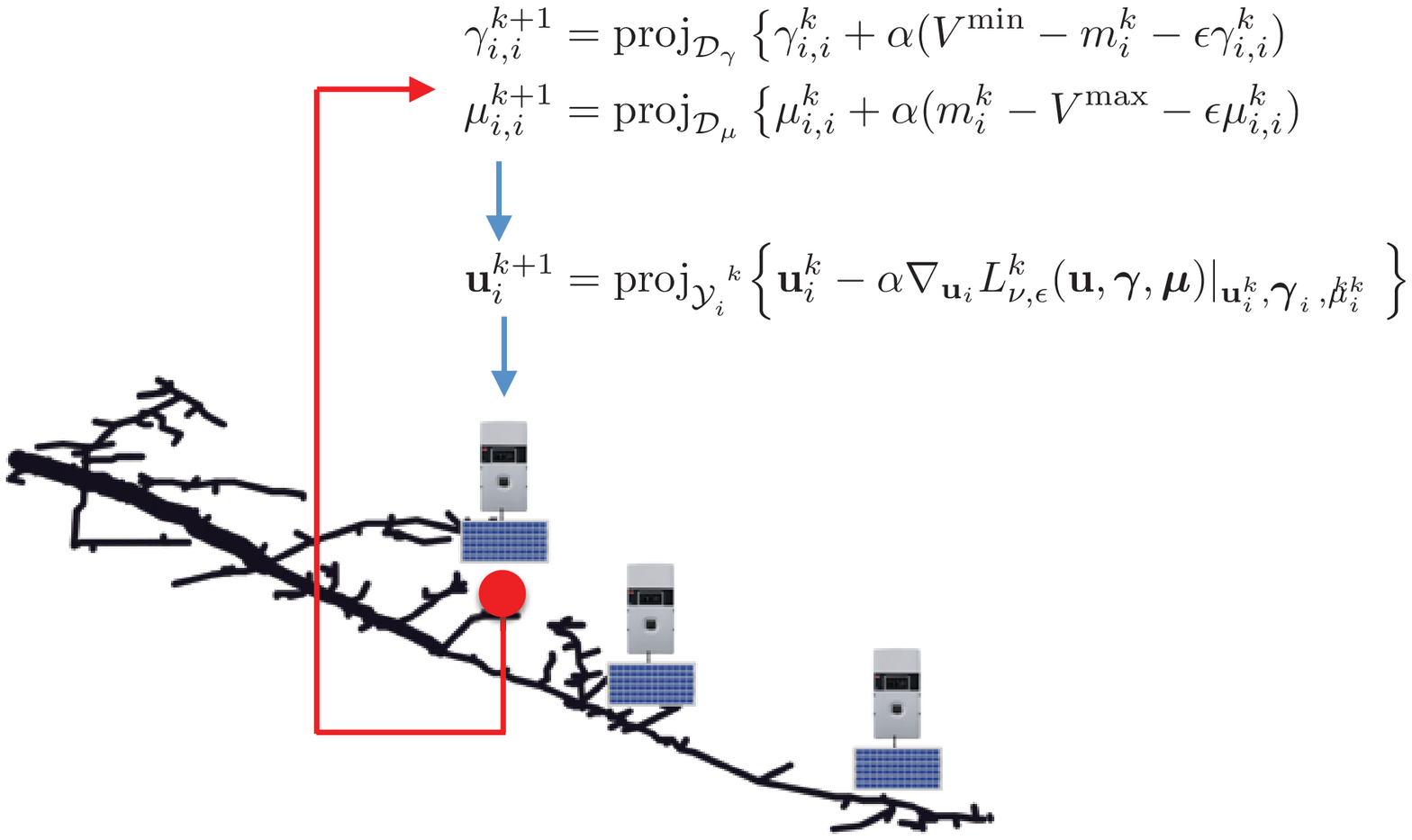} }\vspace{-.2cm}
\caption{(a) Step of the proposed controller architecture, based on voltage measurements gathered at nodes $\cM$. (b) Fast local update performed in between reception of voltage measurements from other nodes of the system. }
\label{F_controller}
\vspace{-.6cm}
\end{figure}

To include actionable feedback from the system, the idea is to replace the algorithmic quantities $\{g^k_n(\bu^k)\}_{n \in \cM}$ and $\{\bar{g}^k_n(\bu^k)\}_{n \in \cM}$ with actual \emph{voltage measurements}; to this end, let $m_n^k$ denote a measurement of the voltage magnitude $|V_n^k|$ acquired at time $k$ from node $n \in \cM$. Further, to account for communication errors in collecting the dual variables at each DER, let $\tilde{\bgamma}_i^{\ell_i(k)}, \tilde{\bmu}_i^{\ell_i(k)}$ represent copies of the most recent multipliers available at DER $i$, with $\ell_i(k) \in \mathbb{N}$ the index of the most recent successful communication. Accordingly, the proposed control architecture amounts to the following iterative steps.

\vspace{.1cm}

\begin{subequations}
\label{eq:controllers}
\noindent \textbf{[S1]} Update power setpoints at each DER $i \in \cG$ as:   
\begin{align} 
& \hspace{-.2cm} \bu_i^{k + 1} \hspace{-.1cm} =  \mathrm{proj}_{\cY_i^{k}} \hspace{-.1cm} \left\{\bu_i^{k}  - \alpha \nabla_{\bu_i}  \cL_{\nu, \epsilon}^k(\bu, \bgamma, \bmu) |_{\bu_i^{k}, \tilde{\bgamma}_i^{\ell_i(k)}, \tilde{\bmu}_i^{\ell_i(k)}} \hspace{-.1cm} \right\}  . \hspace{-.2cm} \label{eq:updatevC} 
\end{align}

\noindent \textbf{[S2]} An aggregator collects voltage measurements $\{m_n^k\}_{n\in\cM}$, updates dual variables as:
\begin{align} 
\gamma_n^{k + 1} & =   \mathrm{proj}_{\cD_\gamma} \left\{\gamma_n^{k} +  \alpha (V^{\mathrm{min}} -  m_n^k - \epsilon \gamma_n^{k}) \right\}  \label{eq:updategammaC} \\
\mu_n^{k + 1} & =   \mathrm{proj}_{\cD_\mu} \left\{\mu_n^{k} +  \alpha (m_n^k - V^{\mathrm{max}}  - \epsilon \mu_n^{k}) \right\}\   \label{eq:updatemuC}
\end{align}
for all $n \in \cM$, and broadcasts dual variables to DERs. 
\end{subequations}

\noindent \textbf{[S3]} Each DER $i \in \cG$ sets the local copies of the dual variables to $\tilde{\bgamma}_i^{k} = \bgamma^{k}$, $\tilde{\bmu}_i^{k} = \bmu^{k}$ if dual variables are received, and $\tilde{\bgamma}_i^{k} = \bgamma_i^{k-1}$,  $\tilde{\bmu}_i^{k} = \bmu_i^{k-1}$ otherwise.

\noindent Go to \textbf{[S1]}. 

\vspace{.1cm}

Steps \textbf{[S1]}--\textbf{[S3]} are illustrated in Figure~\ref{F_controller}(a). Is it worth pointing out that, differently from traditional OPF schemes,~\eqref{eq:controllers} \emph{does not} require knowledge of the loads at locations $\cN \backslash \cG$.  The only information required by the controllers pertains to the line and feeder models, which are utilized to build the  matrices in~\eqref{eq:approximateVoltages}. In the following, the convergence properties of~\eqref{eq:controllers} are analyzed; to this end, pertinent definitions and assumptions are introduced next. 

Let $\bxi_i^k := [\check{\br}_i^k,  \check{\bb}_i^k]^\sfT$, and notice that $\|\xi_i^k \|_2 \leq X_i$ for all $k \geq 0$~\cite{sairaj2015linear,swaroop2015linear}. Further, let $\be_\gamma^k \in \mathbb{R}^M$ and $\be_\mu^k \in \mathbb{R}^M$  collect the dual gradient \emph{errors} $V^{\mathrm{min}} -  y_n^k - \epsilon \gamma_n^{k} - \nabla_{\gamma_n}  \cL_{\nu, \epsilon}^k $ and $y_n^k - V^{\mathrm{max}}  - \epsilon \mu_n^{k} - \nabla_{\mu_n}  \cL_{\nu, \epsilon}^k$, respectively, when actual voltage measurements are utilized instead of the true gradient of the regularized Lagrangian with respect to the dual variables. Then, the following practical assumptions are made.
\vspace{.1cm}

\noindent \emph{Assumption~5}. There exist a constant $e_d \geq 0$ such that $\max\{\|\be_\gamma^k\|_2, \|\be_\mu^k\|_2\}  \leq e_d$ for all $k \geq 0$. \hfill $\Box$ 

\vspace{.1cm}

\noindent \emph{Assumption~6}. For DER $i$, at most $E_i < + \infty$ consecutive communication packets are lost; that is, $\max\{k - \ell_i(k)\} \leq E_i$ for all $k$. \hfill $\Box$ 

\vspace{.1cm}

Under current modeling assumptions, it can be shown that the update~\eqref{eq:updatevC} involves an inexact gradient step, as substantiated in the next lemma.

\vspace{.1cm}

\begin{lemma}
\label{lemma.error}
When $E_i > 0$, one has that $\nabla_{\bu_i}  \cL_{\nu, \epsilon}^k(\bu, \bgamma, \bmu) |_{\bu_i^{k}, \tilde{\bgamma}_i^{\ell_i(k)}, \tilde{\bmu}_i^{\ell_i(k)}}$ is an inexact gradient of the regularized Lagrangian $ \cL_{\nu, \epsilon}^k(\bu_i, \bgamma, \bmu)$ with respect to $\bu_i$ evaluated at $\{\bu_i^{k}, \bgamma^{k}, \bmu^{k}\}$, i.e., $\nabla_{\bu_i}  \cL_{\nu, \epsilon}^k(\bu_i, \bgamma, \bmu) |_{\bu_i^{k}, \tilde{\bgamma}_i^{\ell_i(k)}, \tilde{\bmu}_i^{\ell_i(k)}} = \nabla_{\bu_i}  \cL_{\nu, \epsilon}^k(\bu_i, \bgamma, \bmu) |_{\bu_i^{k}, \bgamma^{k}, \bmu^{k}} + \be^k_{u,i}$, with error bounded as: 
\begin{align} 
\label{eq:error_u} 
\|\be^k_{u,i}\|_2 \leq \alpha E_i X_i [ K + \bar{K} + \epsilon (D_\gamma +  D_\mu) + 2 e_d ] \, .
\end{align}
\vspace{-.2cm}
\hfill $\Box$
\end{lemma}

\vspace{.1cm}

It follows that the overall error in the primal iterate $\be^k_{u} := [(\be^k_{u,1})^\sfT, \ldots, \be^k_{u,N_\cG}]^\sfT$ is bounded too; particularly, \begin{align} 
\label{eq:error_u_total} 
\|\be^k_{u}\|_2 \leq \alpha \left[ \sum_{i \in \cG} E_i^2 X_i^2 [ K + \bar{K} + \epsilon (D_\gamma +  D_\mu) + 2 e_d ]^2 \right]^{\frac{1}{2}} .
\end{align}
Henceforth, denote as $e_u$ the right-hand-side of~\eqref{eq:error_u_total},  and notice that $e_u > e_d$ whenever $E_i > 0$ for all $i \in \cG$. Convergence and tracking properties of the feedback controllers~\eqref{eq:controllers} are established next.

\vspace{.1cm}

\begin{theorem}
\label{theorem.inexact}
Consider the sequence $\{\bz^k\}: = \{\bu^k, \bgamma^k, \bmu^k\}$ generated by~\eqref{eq:controllers}. Let \emph{Assumptions~1--6} hold. For fixed  positive scalars $\epsilon, \nu >0$, if the stepsize $\alpha>0$ is chosen such that
\begin{equation}\label{eq.alpha}
\rho(\alpha):= \sqrt{1 - 2 \eta  \alpha + \alpha^2 L_{\nu, \epsilon}^2} < 1,
\end{equation}
that is $0 <\alpha < 2 \eta/ L_{\nu, \epsilon}^2$, then the sequence $\{\bz^k\}$ converges Q-linearly to $\bz^{*,k} := \{\bu^{*,k}, \bgamma^{*,k}, \bmu^{*,k}\}$ up to the asymptotic error bound given by: 
\begin{align}
\limsup_{k\to \infty} \|\bz^k - \bz^{*,k}\|_2 = \frac{1}{1 - \rho(\alpha)} \Big[\alpha e + \sigma_{\bz}\Big]
\label{eq.asympt_error}
\end{align}
where $e = \sqrt{e_u^2 + 2 e_d^2}$.   \hfill $\Box$
\end{theorem}

Bound~\eqref{eq.asympt_error} can be obtained by following steps similar to~\cite[Thm.~1]{OPFpursuit}, and the proof is omitted due to space limitations; key is to show that, in spite of the the error in the primal updates,~\eqref{eq:controllers} preserves the properties of a strongly monotone operator and leads to a contraction mapping for $\|\bz^k - \bz^{*,k}\|_2$ if~\eqref{eq.alpha} is satisfied.  Equation~\eqref{eq.asympt_error} quantifies the maximum discrepancy between the iterates $ \{\bu^k, \bgamma^k, \bmu^k\}$ generated by the proposed controllers and the (time-varying) minimizer of problem~\eqref{eq:saddlepoint}. From~\cite[Lemma~3.2]{Koshal11} and by using the triangle inequality, a bound for the difference between $\bu^k$ and the solution of~\eqref{Pmg2} can be obtained.   

%

\subsection{Fast local updates}
\label{sec:asynchrounous}

A modified version of the control scheme is proposed next to address the case where communication constraints introduce significant delays in the computation of  steps~\eqref{eq:controllers}. Particularly,~\eqref{eq:controllers} is complemented by local updates of the DER setpoints  based on measurements of voltages at the DER points of connection as described in the following. 

\vspace{.1cm}

\begin{subequations}
\label{eq:controllers_withLocal}
\noindent \textbf{[S1$^\prime$]} Update power setpoints at each DER $i \in \cG$ as:   
\begin{align} 
& \hspace{-.3cm} \bu_i^{k + 1} \hspace{-.1cm} =  \mathrm{proj}_{\cY_i^{k}} \hspace{-.1cm} \left\{\bu_i^{k}  - \alpha \nabla_{\bu_i}  \cL_{\nu, \epsilon}^k(\bu, \bgamma, \bmu) |_{\bu_i^{k}, \tilde{\bgamma}_i^{k}, \tilde{\bmu}_i^{k}} \hspace{-.1cm} \right\} \hspace{-.2cm} \label{eq:updatevC} 
\end{align}

\noindent \textbf{[S2$^\prime$]} An aggregator collects voltage measurements $\{m_n^k\}_{n\in\cM}$, updates dual variables as:
\begin{align} 
\gamma_n^{k + 1} & =   \mathrm{proj}_{\cD_\gamma} \left\{\gamma_n^{k} +  \alpha (V^{\mathrm{min}} -  m_n^k - \epsilon \gamma_n^{k}) \right\}  \label{eq:updategammaC} \\
\mu_n^{k + 1} & =   \mathrm{proj}_{\cD_\mu} \left\{\mu_n^{k} +  \alpha (m_n^k - V^{\mathrm{max}}  - \epsilon \mu_n^{k}) \right\}\   \label{eq:updatemuC}
\end{align}
for all $n \in \cM$, and broadcasts dual variables to DERs. 

\noindent \textbf{[S3$^\prime$]} At each DER $i \in \cG$, update the local copies of the dual variables as:   

\noindent $\bullet$ If  $\bgamma^{k}$ and $\bmu^{k}$ are available and are received, set $\tilde{\bgamma}_i^{k} = \bgamma^{k}$, $\tilde{\bmu}_i^{k} = \bmu^{k}$;

\noindent $\bullet$ If  $\tilde{\bgamma}_i^{k}$ and $\tilde{\bmu}^{k}$ are not available, measure the voltage magnitude $|V_i^k|$ at the point of connection and update the $i$th entry of $\tilde{\bgamma}_i^{k}$ and $\tilde{\bmu}_i^{k}$ as
\begin{align} 
\tilde{\gamma}_{i,i}^{k + 1} & =   \mathrm{proj}_{\cD_\gamma} \left\{\tilde{\gamma}_{i,i}^{k} +  \alpha (V^{\mathrm{min}} -  m_i^k - \epsilon \tilde{\gamma}_{i,i}^{k}) \right\}  \label{eq:updategammaClocal} \\
\tilde{\mu}_{i,i}^{k + 1} & =   \mathrm{proj}_{\cD_\mu} \left\{\tilde{\mu}_i^{k} +  \alpha (m_i^k - V^{\mathrm{max}}  - \epsilon \tilde{\mu}_{i,i}^{k}) \right\} .  \label{eq:updatemuClocal}
\end{align}

\noindent The remaining entries are not updated; i,e.,  $\tilde{\gamma}_{i,j}^{k + 1} = \tilde{\gamma}_{i,j}^{k}$ and $\tilde{\mu}_{i,j}^{k + 1} = \tilde{\mu}_{i,j}^{k}$ for all $j \in \cM \backslash \{i\}$. 

\noindent Go to \textbf{[S1$^\prime$]}. 
\end{subequations}

\vspace{.1cm}

As shown in Figure~\ref{F_controller}, steps  \textbf{[S1$^\prime$]}--\textbf{[S3$^\prime$]} allow each DER inverter $i$ to update the setpoints $\bu_i^{k}$ at a \emph{faster} time scale, based on local measurements of the voltage level at the DER point on interconnection; et each time step, DER $i$ continuously updates the $i$th entry of $\tilde{\bgamma}_i^{k}$ and $\tilde{\bmu}_i^{k}$ and  computes the setpoints $\bu_i^{k}$. The remaining entries of $\tilde{\bgamma}_i^{k}$ and $\tilde{\bmu}_i^{k}$ are updated when the  vectors $\bgamma^{k}$ and $\bmu^{k}$ become available.  It is worth emphasizing that  steps  \textbf{[S1$^\prime$]}--\textbf{[S3$^\prime$]} consider the case where an aggregator collects voltage measurements and broadcasts the updated version of the dual variables\verb""; the algorithm can be suitably modified to account for the case where each DER receives measurements of the voltage across nodes $n \in \cM$, and updates the local copies of the dual variables based on $\{m_n^k\}_{n \in \cM}$. This scenario leads to an operational setup where at each time $k$ the DER updates the entries of $\tilde{\bgamma}_i^{k}$ and $\tilde{\bmu}_i^{k}$ that correspond to the subset of nodes from which voltage measurements are received. 

The results of Lemma~\ref{lemma.error} and Theorem~\ref{theorem.inexact} can be adapted to \textbf{[S1$^\prime$]}--\textbf{[S3$^\prime$]}. In this case, $M_i$ represents the number of iterations that are necessary for DER $i$ to update all the entries of the dual variables (or to receive measurements of all voltages in $\cM$).

\begin{figure}[t] 
\centering
\psfrag{1}[c]{\footnotesize 1}
\psfrag{2}[c]{\footnotesize 2}
\psfrag{3}[c]{\footnotesize 3}
\psfrag{4}[c]{\footnotesize 4}
\psfrag{5}[c]{\footnotesize 5}
\psfrag{6}[c]{\footnotesize 6}
\psfrag{7}[c]{\footnotesize 7}
\psfrag{8}[c]{\footnotesize 8}
\psfrag{9}[c]{\footnotesize 9}
\psfrag{10}[c]{\footnotesize 10}
\psfrag{11}[c]{\footnotesize 11}
\psfrag{12}[c]{\footnotesize 12}
\psfrag{13}[c]{\footnotesize 13}
\psfrag{14}[c]{\footnotesize 14}
\psfrag{15}[c]{\footnotesize 15}
\psfrag{16}[c]{\footnotesize 16}
\psfrag{17}[c]{\footnotesize 17}
\psfrag{18}[c]{\footnotesize 18}
\psfrag{19}[c]{\footnotesize 19}
\psfrag{20}[c]{\footnotesize 20}
\psfrag{21}[c]{\footnotesize 21}
\psfrag{22}[c]{\footnotesize 22}
\psfrag{23}[c]{\footnotesize 23}
\psfrag{24}[c]{\footnotesize 24}
\psfrag{25}[c]{\footnotesize 25}
\psfrag{26}[c]{\footnotesize 26}
\psfrag{27}[c]{\footnotesize 27}
\psfrag{28}[c]{\footnotesize 28}
\psfrag{29}[c]{\footnotesize 29}
\psfrag{30}[c]{\footnotesize 30}
\psfrag{31}[c]{\footnotesize 31}
\psfrag{32}[c]{\footnotesize 32}
\psfrag{33}[c]{\footnotesize 33}
\psfrag{34}[c]{\footnotesize 34}
\psfrag{35}[c]{\footnotesize 35}
\psfrag{36}[c]{\footnotesize 36}
\psfrag{37}[c]{\footnotesize 37}
%
\vspace{.25cm}
\includegraphics[width=.35\textwidth]{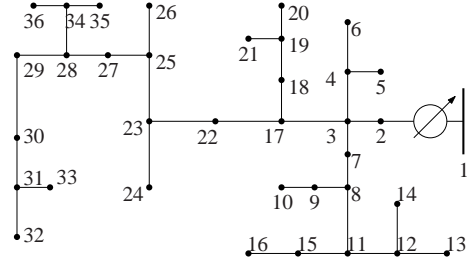}
\caption{IEEE 37-node feeder.}
\label{F_feeder}
\end{figure}

\section{Example of Application}
\label{sec:application}

Consider a modified version of the IEEE 37-node test feeder shown in Figure~\ref{F_feeder}. The modified network is obtained by considering a single-phase equivalent, and by replacing the loads specified in the original dataset with real load data measured from feeders serving a neighborhood called Anatolia in CA during the week of August 2012~\cite{Bank13}. Particularly, the data have a granularity of $1$ second, and represent the loading of secondary transformers. Line impedances, shunt admittances, as well as active and reactive loads are adopted from the original dataset. With reference to Fig.~\ref{F_feeder}, it is assumed that  PV systems are located at nodes $4$, $7$, $10$, $13$, $17$, $20$, $22$, $23$, $26$, $28$, $29$, $30$, $31$, $32$, $33$, $34$, $35$, and $36$, and their generation profile is simulated based on the measured solar irradiance data available in~\cite{Bank13}. Solar irradiance data have a granularity of $1$ second. The rating of the PV inverters are $300$ kVA for $i = 3$, $350$ kVA for $i = 15, 16$, and $200$ kVA for the remaining PV inverters.

\begin{figure}[t] 
\subfigure[]{\hspace{-.2cm}\includegraphics[width=9.0cm]{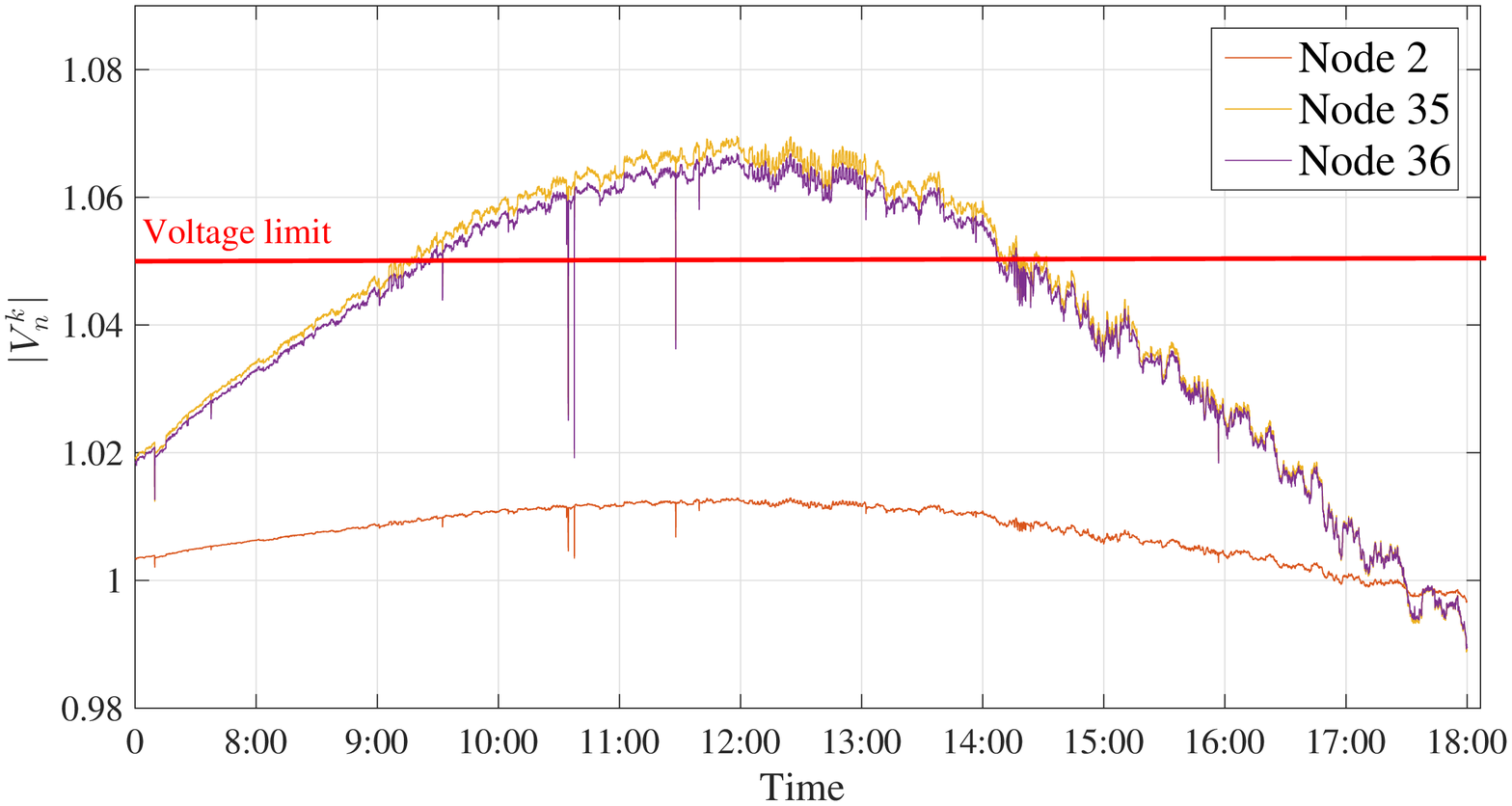} }
\subfigure[]{\hspace{-.2cm}\includegraphics[width=9.0cm]{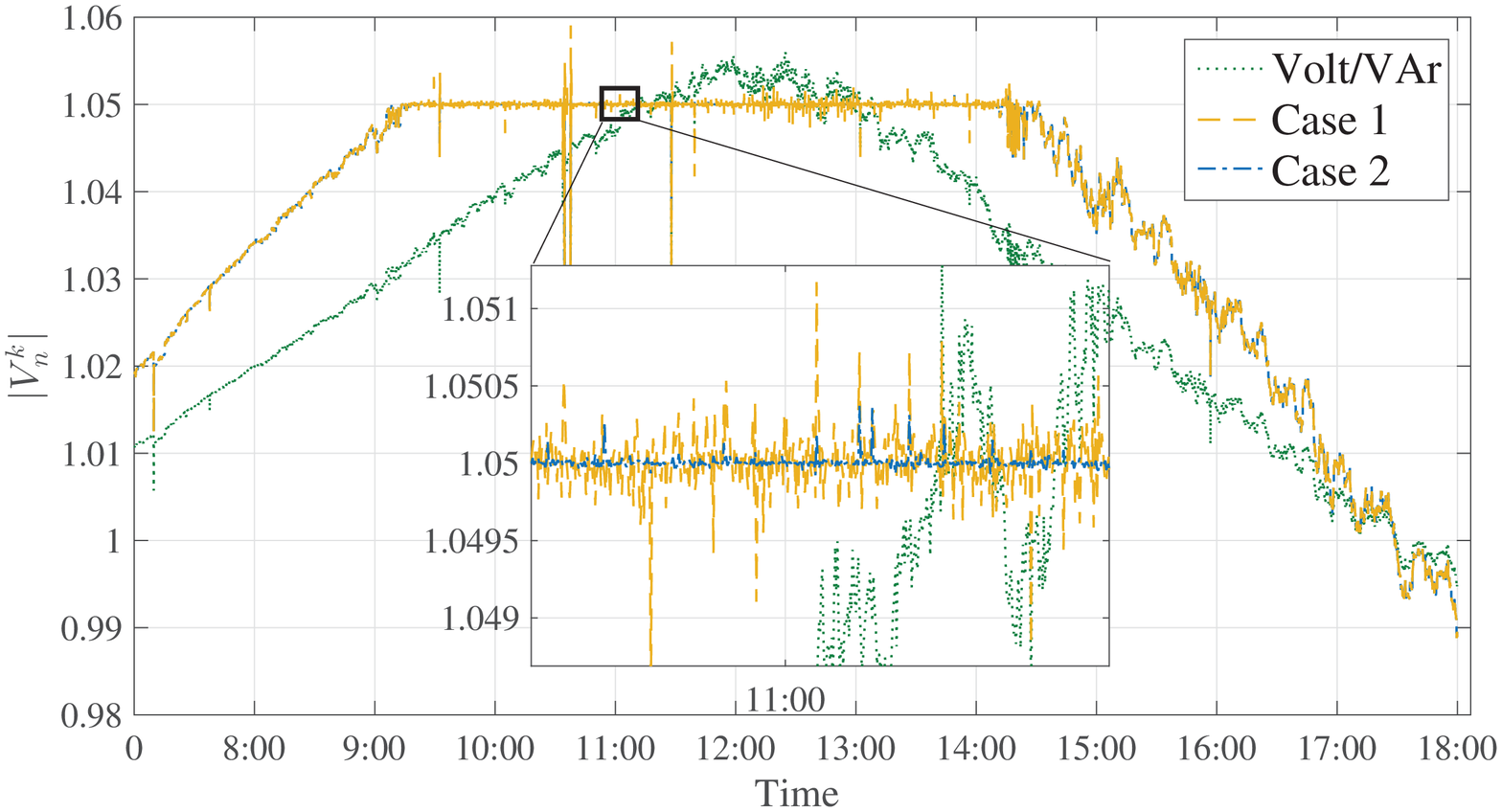} }
\subfigure[]{\hspace{-.2cm}\includegraphics[width=9.0cm]{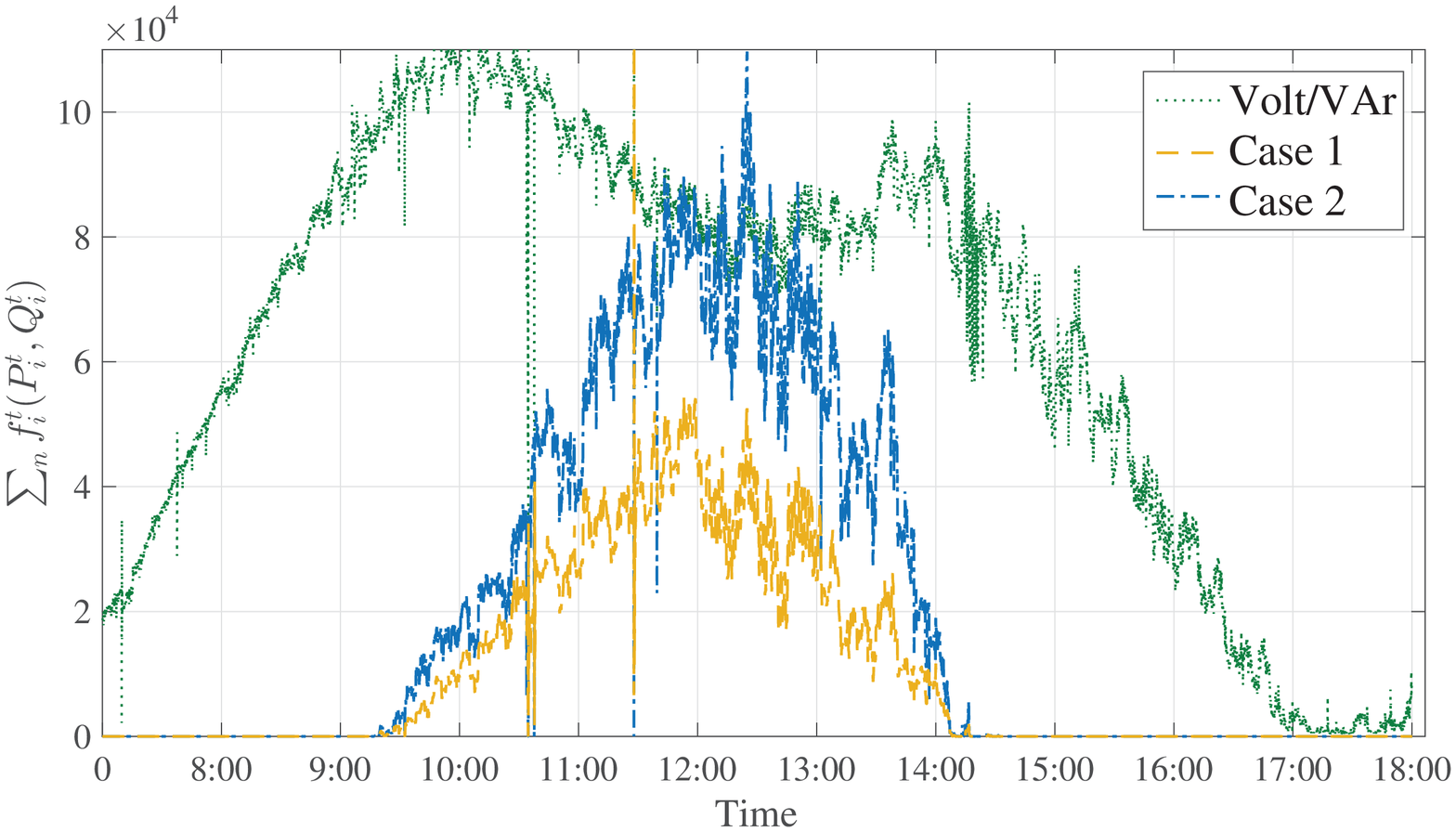} }\vspace{-.2cm}
\caption{(a) Voltage profile at representative nodes when no control is implemented at the PV systems. (b) Achieved voltage profile with Volt/VAr local control and proposed controllers. (c) Cost of ancillary service provisioning.}
\label{F_voltage}
\vspace{-.6cm}
\end{figure}

The goal of this simulation study is to demonstrate how the proposed controllers can reliably prevent overvoltages that are likely to be experienced during periods when PV generation exceeds the demand. The minimum and maximum voltage limits are set to $0.95$ pu and $1.05$ pu, respectively. With this simulation setup, when no actions are taken to prevent overvoltages, one would have voltages well beyond $1.05$ pu at nodes $29$--$36$, with the most severe overvoltage conditions experiences at node $35$. This is clear from the voltage profile provided in Figure~\ref{F_voltage}(a).

Two cases are considered: 

\emph{Case 1}: Controller \textbf{[S1]}--\textbf{[S3]} is implemented, and the primal-dual updates represented in Figure~\ref{F_controller}(a) are repeated every $1$ second. 

\emph{Case 2}: Controller \textbf{[S1$^\prime$]}--\textbf{[S3$^\prime$]}, where the global steps represented in Figure~\ref{F_controller}(a) are repeated every $1$ second,  and are complemented by the local steps~\eqref{eq:updatevC} and~\eqref{eq:updategammaClocal}--\eqref{eq:updatemuClocal}; the steps represented in Figure~\ref{F_controller}(b) are performed every $0.1$ seconds. This way, $E_i$ turns out to be $E_i = 9$ for all $i \in \cG$. 

The target optimization objective~\eqref{mg-cost2} is set to $f_n^k(\bu_n^k) =  c_q (Q_n^k)^2 + c_p (P_{\textrm{av},n}^k - P_n^k)^2$ to minimize the amount of real power curtailed from the PV systems and to minimize the amount of reactive power injected or absorbed. The coefficients are set to $c_p = 3$ and $c_q = 1$ for all PV systems to discourage real power curtailment. It is assumed that the dual ascent step is performed at the utility/aggregator, which subsequently broadcasts the dual variables to the PV systems. The  controller parameters are set as $\nu = 10^{-3}$, $\epsilon = 10^{-4}$, and $\alpha = 0.2$. The stepsize $\alpha$ was selected experimentally.  

In \emph{Case 2}, the PV system setpoints are updated at a faster time scale by utilizing local voltage measurements [cf. Figure~\ref{F_controller}(b)]; voltage across the network are collected every $1$ s, and are utilized to update  all the entries of the dual variables [cf. Figure~\ref{F_controller}(a)]. The performance of the proposed controllers is compared with local Volt/VAr control.This involves a linear trade off, where inverters set $Q_n^k = 0$ when $|V_n^k| = 1$ pu and linearly increase the reactive power to $Q_n^k = -\sqrt{S_n^2 - (P_{\textrm{av},n}^t)^2}$ when $|V_n^k| \geq 1.05$ pu. The PV-inverters measure the voltage magnitude and update the reactive-power setpoint every 0.1 seconds. 

Figure~\ref{F_voltage}(b) illustrates the voltage profile obtained at node 35 using the proposed controllers as well as local Volt/VAr control. First, it can be seen that Volt/VAr control fails in resolving overvoltage conditions in the considered setup. In contrast, the proposed controllers ensure that voltage limits are satisfied. However, it can be clearly seen that the controllers in \emph{Case 2} yield  a smoother voltage profile, and this ensures higher power-quality  at both the customer and utility sides. 

Figure~\ref{F_voltage}(c) reports the cost achieved by the proposed controllers; that is $\sum_{n \in \cG} c_q (Q_n^k)^2 + c_p (P_{\textrm{av},n}^k - P_n^k)^2$, for all $k$. This is compared with the cost of reactive power provisioning incurred by Volt/VAr control, which is computed as $\sum_{n \in \cG} c_q (Q_n^k)^2$, for all $k$. The advantages of the proposed controllers are evident, as they enable voltage regulation at a lower cost. However, it can be seen that the improved voltage profile obtained in \emph{Case 2} comes at a higher cost.

\section{Concluding Remarks}
\label{sec:conclusions}

This paper addressed the synthesis of feedback controllers that seek DER setpoints corresponding to AC OPF solutions. Appropriate linear approximations of the AC power flow equations were utilized to facilitate the development of low-complexity controllers; and primal-dual methods were leveraged for the controller synthesis.  The tracking capabilities of the proposed controllers were analytically established and numerically corroborated for the case of communication-packet losses and partial updates of control signals. 

\appendix

\emph{Proof of Lemma}~\ref{lemma.error}. Define $\bzeta_i^{k,k-i} := \nabla_{\bu_i}  \cL_{\nu, \epsilon}^k(\bu, \bgamma, \bmu) |_{\bu_i^{k}, \bgamma^{k-i}, \bmu^{k-i}}$, which is given by
\begin{align} 
\bzeta_i^{k,k-i} \hspace{-.1cm} & = \hspace{-.1cm} \nabla_{\bu_i}[f_1^k(\bu_1), \ldots, f_{N_\cG}^k(\bu_{N_\cG})]^\sfT |_{\bu_i^{k}} \nonumber \\
& \hspace{1.8cm} + \bxi_i^k (\bmu^{k-i}  - \bgamma^{k-i}) + \nu \bu_i^{k} .
\end{align}
Recall that $\|\bgamma_i^k\|_2 \leq D_\gamma$ and $\|\bmu_i^k\|_2 \leq D_\mu$ for all $k \geq 0$, and notice that the norm of the vector $\bxi_i^k = [\check{\br}_i^k,\check{\bb}_i^k]^\sfT$ can be bounded as $\|\xi_i^k \|_2 \leq X_i$ for all $k \geq 0$~\cite{sairaj2015linear}. Next, notice that $\bzeta_i^{k,k-M}$ can be written as
\begin{align} 
\label{eq:p_l_1}
\bzeta_i^{k,k-M}  = \bzeta_i^{k,k} + \be^k_{u,i} \, ,
\end{align}
where $ \be^k_{u,i}  = \sum_{j = 1}^{E_i} (\bzeta_i^{k,k-j} - \bzeta_i^{k,k-j+1})$. Expanding on~\eqref{eq:p_l_1}, one obtains that $\be^k_{u,i} = \bxi_i^k \sum_{j = 1}^{E_i} [(\bmu_i^{k-j} - \bmu_i^{k-j+1}) + (\bgamma^{k-j+1} - \bgamma^{k-j})]$. Then, using the triangle inequality, one has that 
$\|\be^k_{u,i} \|_2 \leq X_i \sum_{j = 1}^{E_i} [\|\bmu^{k-j} - \bmu^{k-j+1}\|_2 + \|\bgamma^{k-j+1} - \bgamma^{k-j}\|_2]$ .
Next, $\|\bmu^{k-j} - \bmu^{k-j+1}\|_2$ can be bounded as:
\begin{subequations}
\begin{align} 
& \|\bmu^{k-i} - \bmu^{k-i+1}\|_2 \nonumber \\
& = \|\bmu^{k-i} \hspace{-.1cm} - \hspace{-.1cm}  \mathrm{proj}_{\cD_\mu} \hspace{-.1cm}  \left\{\bmu_n^{k-i} +  \alpha ( \bm^k - \mathbf{1} V^{\textrm{max}} - \epsilon \bmu^{k-i}) \hspace{-.1cm}\right\} \|_2  \hspace{-.1cm}  \label{eq:p_l_3} \\
& = \|\bmu^{k-i} \hspace{-.1cm} - \hspace{-.1cm}  \mathrm{proj}_{\cD_\mu} \hspace{-.1cm}  \left\{\bmu_n^{k-i} +  \alpha ( \bar{\bg}^{k-i}(\bu^{k})  + \be_\mu^{k-i} - \epsilon \bmu^{k-i}) \right\} \|_2  \hspace{-.1cm}  \label{eq:p_l_4} \\
& \leq \|  \alpha ( \bar{\bg}^{k-i}(\bu^{k})  + \be_\mu^{k-i} - \epsilon \bmu^{k-i} )  \|_2 \label{eq:p_l_5} \\
& \leq   \alpha ( \|\bar{\bg}^{k-i}(\bu^{k})\|_2 + \|\be_\mu^{k-i}\|_2 + \epsilon \| \bmu^{k-i}\|_2 )  \label{eq:p_l_6}  \\
& \leq \alpha (\bar{K} + e_d + \epsilon D_\mu ) \,  \label{eq:p_l_7} 
\end{align}
\end{subequations}
where $\bm^k$ in~\eqref{eq:p_l_3}  collects all the voltage measurements $\bm^k := [m_1^k, \ldots, m_{M}^k]^\sfT$, and the non-expansive property of the projection operator, along with the fact that $\bmu^{j} = \mathrm{proj}_{\cD_\mu}\{\bmu^{j}\}$, is utilized to derive~\eqref{eq:p_l_5}.  Using~\eqref{eq:p_l_5}, it follows that $\sum_{j = 1}^{E_i} [\|\bmu^{k-j} - \bmu^{k-j+1}\|_2 \leq \alpha {E_i} (\bar{K} + e_d + \epsilon D_\mu )$. Following similar steps, one can show that $\sum_{j = 1}^{M_i} [\|\bgamma^{k-j} - \bgamma^{k-j+1}\|_2 \leq \alpha E_i (K + e_d + \epsilon D_\gamma )$. Bound~\eqref{eq:error_u} readily follows.

\bibliographystyle{IEEEtran}
\bibliography{biblio.bib}

\end{document}